\documentclass[a4paper,10pt]{amsart}
\usepackage{amssymb,amsmath,amsfonts,amsthm}
\usepackage[dvips]{graphicx}
\usepackage{psfrag}

\theoremstyle{plain}
\newtheorem{main}{Theorem}

\newtheorem{theorem}{Theorem}[section]
\newtheorem{lemma}[theorem]{Lemma}
\newtheorem{proposition}[theorem]{Proposition}
\newtheorem{corollary}[theorem]{Corollary}
\theoremstyle{remark}
\newtheorem{remark}[theorem]{Remark}

\newtheorem{example}[theorem]{Example}

\newcommand{\quand}{\quad\text{and}\quad}
\newcommand{\per}{\operatorname{per}}

\def\SL{\operatorname{SL}}
\def\supp{\operatorname{supp}}

\def\cM{{\mathcal M}}

\def\cR{{\mathcal R}}

\def\cW{{\mathcal W}}

\def\id{\operatorname{id}}

\def\loc{{\operatorname{loc}}}

\def\RR{{\mathbb R}}

\def\ZZ{{\mathbb Z}}
\def\NN{{\mathbb N}}

\def\hA{\hat A}
\def\hB{\hat B}
\def\hE{\hat E}
\def\hF{\hat F}

\def\hM{{\hat M}}
\def\hV{\hat V}

\def\tW{\widetilde W}
\def\hf{\hat f}

\def\hm{{\hat m}}
\def\hp{{\hat p}}
\def\hq{{\hat q}}
\def\hx{{\hat x}}
\def\hy{{\hat y}}
\def\hz{{\hat z}}
\def\hxi{{\hat\xi}}

\def\hmu{{\hat\mu}}

\def\pR{\mathbb P\mathbb R}

\def\H{{1+\epsilon}}

\title{Continuity of Lyapunov exponents\\ in the $C^0$ topology}
\author{Marcelo Viana and Jiagang Yang}
\date{\today}
\thanks{M.V. and J.Y. were partially supported by CNPq, FAPERJ, and PRONEX}

\address{IMPA, Est. D. Castorina 110, 22460-320 Rio de Janeiro, Brazil}
\email{viana\@@impa.br}

\address{Departamento de Geometria, Instituto de Matem\'atica e Estat\'\i stica, Universidade Federal Fluminense, Niter\'oi, Brazil}
\email{yangjg\@@impa.br}

\begin{document}

\begin{abstract}
We prove that the Bochi-Ma\~{n}\'{e} theorem is false, in general, for linear cocycles
over non-invertible maps: there are $C^0$-open subsets of linear cocycles that are not
uniformly hyperbolic and yet have Lyapunov exponents bounded from zero.
\end{abstract}

\maketitle

\hfill Dedicated to Welington de Melo

\setcounter{tocdepth}{1} \tableofcontents

\section{Introduction}\label{s.introduction}

Bochi~\cite{Bochi-tese,Boc02} proved that every continuous $\SL(2)$-cocycle over an aperiodic
invertible system can be approximated in the $C^0$ topology by cocycles whose Lyapunov exponents
vanish, unless it is uniformly hyperbolic. The (harder) version of the statement for
derivative cocycles of area-preserving diffeomorphisms on surfaces had been claimed by
Ma\~{n}\'{e}~\cite{Man83} almost two decades before.
Bochi~\cite{Bochi-tese,Boc02} also completed the proof of this harder claim, based on an outline by
Ma\~n\'{e}. These results were then extended to arbitrary dimension by Bochi, Viana~\cite{BcV05}
and Bochi~\cite{Boc09}.

In this paper, we prove that the Bochi-Ma\~{n}\'{e} theorem does not hold, in general,
for cocycles over non-invertible systems: surprisingly, in the non-invertible setting
there exist $C^0$-open sets of $\SL(2)$-cocycles whose exponents are bounded away from zero.
Indeed, we provide two different constructions of such open sets.

The first one (Theorem~\ref{t.main1}) applies to H\"{o}lder continuous cocycles satisfying
a bunching condition. The second one (Theorem~\ref{t.main2}) has no bunching hypothesis but
requires the cocycle to be $C^\H$ and to be hyperbolic at some periodic point.
In either case, we assume some form of irreducibility. A suitable extension of the Invariance
Principle (Bonatti, Gomez-Mont, Viana~\cite{BGV03}, Avila, Viana~\cite{Extremal}) that we prove
here gives that these cocycles have non-zero Lyapunov exponents.
We also prove that they are continuity points for the Lyapunov exponents, relative to the
$C^0$ topology, and thus the claim follows.

Continuity of Lyapunov exponents, especially with respective to finer topologies, has been the
object of considerable recent interest.  See Viana~\cite[Chapter~10]{LLE},
Duarte, Klein~\cite{DuK16b} and references therein. It was conjectured by Viana~\cite{LLE}
that Lyapunov exponents are always continuous among H\"older continuous \emph{fiber-bunched}
$\SL(2)$-cocycles, and this has just been proved by Backes, Brown, Butler~\cite{BBB}.
In fact, they prove a stronger conjecture also from Viana~\cite{LLE}: Lyapunov exponents vary
continuously on any family of $\SL(2)$-cocycles with continuous invariant holonomies.
Improving on a construction of Bocker, Viana in \cite[Chapter~9]{LLE}, Butler~\cite{But} also
shows that the fiber-bunching condition is sharp for continuity in some cases.

These and many other related results require the cocycles to some fair amount of regularity,
starting from H\"older continuity. In view also of the Bochi-Ma\~n\'e theorem, continuity in the
$C^0$ topology (outside the uniformly hyperbolic realm) as we exhibit here, comes as a bit
of a surprise.
It seems that the explanation may lie on the fact that existence of an invariant stable holonomy
comes for free in the non-invertible case. 

\section{Statement of results}\label{s.statement}

Let $f:M\to M$ be a continuous uniformly expanding map on a compact metric space.
By this we mean that there are $\rho>0$ and $\sigma>1$ such that, for any $x\in M$,
\begin{enumerate}
\item[(i)] $d(f(x),f(y)) \ge\sigma d(x,y)$ for every $y \in B(x,\rho)$ and
\item[(ii)] $f(B(x,\rho))$ contains the closure of $B(f(x),\rho)$.
\end{enumerate}
%
Take $f$ to be topologically mixing and let $\mu$ be the equilibrium state of some H\"{o}lder
continuous potential (see \cite[Chapter 11]{FET}). Then $\mu$ is $f$-invariant and ergodic,
and the support is the whole $M$.

Let $\hf:\hM\to\hM$ be the natural extension of $f$, that is, the shift map
$$
(\dots, x_{-n}, \dots, x_{-1}, x_0) \mapsto (\dots, x_{-n}, \dots, x_{-1}, x_0, f(x_0))
$$
in the space $\hM$ of sequences $(x_{-n})_{n}$ such that $f(x_{-n})=x_{-n+1}$ for every $n\ge 1$.
Then $\hf$ is a hyperbolic homeomorphism (see \cite[Definition~1.3 and Section~6]{Almost}).
For any $\hx=(x_{-n})_n$ in $\hM$,
\begin{itemize}
\item the local stable set $W_\loc^s(\hx)$ is the fiber $\pi^{-1}(\hx)$ of the
canonical projection $\pi(\hx)=x_0$, and
\item the local unstable set $W_\loc^u(\hx)$ consists of the points $\hy=(y_{-n})_n$ such that
$d(x_{-n},y_{-n})<\rho$ for every $n\ge 0$.
\end{itemize}

Let $\hmu$ be the lift of $\mu$ to $\hM$, that is, the unique $\hf$-invariant measure
that projects down to $\mu$ under $\pi$.
Then $\hmu$ is ergodic and supported on the whole $\hM$, and it has local product
structure (see ~\cite[Section~2.2]{BoV04}).

The \emph{projective cocycle} defined by a continuous  map $A:M\to\SL(2)$ over the
transformation $f$ is the map $F_A:M\times\pR^2\to M\times\pR^2$,
$F_A(x,v) = (f(x), A(x)v)$. Denote $A^n(x)=A(f^{n-1}(x))\cdots A(x)$ for every $n\ge 1$.
By~\cite{FK60,Ki68}, there exists $\lambda(A) \ge 0$, called \emph{Lyapunov exponent},
such that
\begin{equation}\label{eq.exponents}
\lim_n \frac 1n \log \|A^n(x)\|
= \lim_n \frac 1n \log \|A^n(x)^{-1}\|
= \lambda(A)
\text{ for $\mu$-almost every $x\in M$.}
\end{equation}

We say that $A$ is \emph{$u$-bunched} if there exists $\theta\in(0,1]$ such that $A$ is
$\theta$-H\"{o}lder continuous and
\begin{equation}\label{eq.ubunched}
\|A(x)\| \|A(x)^{-1}\| \, \sigma^{-\theta} < 1
\text{ for every } x \in M.
\end{equation}
See \cite[Definitions~1.11 and 2.2]{BGV03} and \cite[Definition~2.2 and Remark~2.3]{ASV13}.

\begin{remark}\label{r.holder}
The function $d_\theta(x,y)=d(x,y)^\theta$ is also a distance in $M$, and it satisfies (i) above with
$\sigma$ replaced with $\sigma^\theta$. Moreover, $A$ is $1$-H\"{o}lder continuous with respect
to $d_\theta$ if it is $\theta$-H\"{o}lder continuous with respect to $d$.
Thus, up to replacing the distance in $M$, it is no restriction to suppose that $\theta=1$,
and we do so.
\end{remark}

Let $\hA:\hM \to \SL(2)$ be defined by $\hA = A \circ \pi$. Assuming that $A$
is $u$-bunched, the cocycle $\hF_A$ defined by $\hA$ over $\hf$ admits invariant
$u$-holonomies (see \cite[Section~1.4]{BGV03} and \cite[Section~3]{ASV13}), namely,
$$
h_{\hx,\hy}^u = \lim_n \hA^n\big(\hf^{-n}(\hy)\big) \hA^n\big(\hf^{-n}(\hx)\big)^{-1}
\text{ for any $\hy\in W^u_\loc(\hx)$.}
$$
As $\hA$ is constant on local stable sets, $\hF_A$ also admits trivial
invariant $s$-holonomies:
$$
h_{\hx,\hy}^s = \id
\text{ for any $\hy\in W^s_\loc(\hx)$.}
$$

\begin{remark}\label{r.sbunched}
It is not difficult to find a distance on $\hM$ relative to which $\hA$ is $s$-bunched,
in addition to being $u$-bunched.
\end{remark}

A probability measure $\hm$ on $\hM\times\pR^2$ is said to be \emph{$u$-invariant} if it admits
a disintegration $\{\hm_\hx: \hx \in\hM\}$ along the fibers $\{\hx\}\times\pR^2$ such that
\begin{equation}\label{eq.u-invariant}
\big(h_{\hx,\hy}^u\big)_*\hm_{\hx} = \hm_{\hy}
\text{ for any $\hy \in W^u_\loc(\hx)$.}
\end{equation}
Similarly, we say that $\hm$ is \emph{$s$-invariant} if it admits a disintegration $\{\hm_\hx: \hx \in\hM\}$
along the fibers $\{\hx\}\times\pR^2$ such that
\begin{equation}\label{eq.s-invariant}
\hm_{\hx} = \hm_{\hy}
\text{ for any $\hy \in W^s_\loc(\hx)$.}
\end{equation}
A $u$-invariant (respectively $s$-invariant) probability measure $\hm$ is called a \emph{$u$-state}
(respectively, an $s$-state) if it is also invariant under $\hF_A$.
We call $\hm$ an \emph{$su$-state} (see \cite[Section~4]{Extremal}) if it is both a $u$-state and
an $s$-state.

\begin{main}\label{t.main1}
If $A$ is $u$-bunched and has no $su$-states then $\lambda(A) >0$ and $A$ is a continuity point for
the function $B \mapsto \lambda(B)$ in the space of continuous maps $B: M\to\SL(2)$ equipped with
the $C^0$ topology. In particular, the Lyapunov exponent $\lambda(\cdot)$ is bounded away from zero
on a $C^0$-neighborhood of $A$.
\end{main}

\begin{example}\label{ex.no_su-states}
Take $f:S^1\to S^1$, $f(x) = k x \mod \ZZ$, for some integer $k \ge 2$,
and $\mu$ to be the Lebesgue measure on $S^1$. Let $A:S^1\to\SL(2)$ be given by
$A(x) = A_0 R_x$, where $A_0 \in \SL(2)$ is a hyperbolic matrix and $R_x$ is the rotation by
angle $x$. $A$ is $1$-H\"{o}lder continuous and, in view of the definition \eqref{eq.ubunched},
it is $u$-bunched provided $k > \|A_0\| \|A_0^{-1}\|$.

We claim that $A$ has no $su$-states if $k$ is large enough. Indeed, suppose that $\hF_A$ has some
$su$-state $\hm$. Then, by \cite[Proposition~4.8]{Extremal}, $\hm$ admits a continuous disintegration
$\{\hm_\hx:\hx\in\hM\}$ which is simultaneously $s$-invariant, $u$-invariant and $\hF_A$-invariant.
By $s$-invariance, we may write the disintegration as $\{\hm_x:x\in M\}$ instead.
Continuity and invariance under the dynamics imply that $(A_0)_* \hm_0 = \hm_0$.
Since $A_0$ is hyperbolic, this means that $m_0$ is either a Dirac mass or a convex combination
of two Dirac masses. Thus, by holonomy invariance, either every $\hm_x$ is a Dirac mass or every
$\hm_x$ is supported on exactly $2$ points.

In the first case, $\xi(x) = \supp\hm_x$ defines a map $\xi:S^1\to\pR^2$ which is continuous
and invariant under the cocycle:
\begin{equation*}
\xi(f(x)) = A_0 R_x \xi(x) \text{ for every $x \in S^1$.}
\end{equation*}
It follows that the degree $\deg\xi$ satisfies $k \deg \xi = \deg \xi + 2$ (the term $2$ comes from the fact
that $S^1\to\pR^2$, $x \mapsto R_x v$ has degree $2$ for any $v$).
This is impossible when $k\ge 4$, and so this first case can be disposed of.
In the second case, $\xi(x) = \supp\hm_x$ defines a continuous invariant section with values in
the space $\pR^{2,2}$ of pairs of distinct points in $\pR^2$.
This can be reduced to the previous case by considering the $2$-to-$1$ covering map $S^1 \to S^1$,
$z\mapsto 2z \mod\ZZ$ (notice that $f$ is its own lift under this covering map).
Thus, this second case can also be disposed of. This proves our claim that $A$ has no $su$-states.

By Theorem~\ref{t.main1}, it follows that $\lambda(B)>0$ for every continuous $B:S^1\to\SL(2)$
in a $C^0$-neighborhood of $A$. Incidentally, this shows that \cite[Corollary~12.34]{Beyond} is not
correct: indeed, the ``proof'' assumes the Bochi-Ma\~n\'{e} theorem in the non-invertible case.
\end{example}

Now let $f:M\to M$ be a $C^\H$ (that is, $C^1$ with H\"{o}lder continuous derivative)
expanding map on a compact manifold $M$ and $A:M\to\SL(2)$ be a $C^\H$ function.
All the other objects, $\mu$, $F_A$, $\pi$, $\hM$, $\hf$, $\hmu$, $\pi$,
$\hA$ and $\hF_A$  are as before. An \emph{invariant section} is a continuous map
$\hxi:\hM\to\pR^2$ or $\hxi:\hM\to\pR^{2,2}$ such that
$$
\hA(\hx)\hxi(\hx)=\hxi(\hf(\hx)) \text{ for every $\hx\in\hM$.}
$$

\begin{main}\label{t.main2}
If $A$ admits no invariant section then it is continuity point for the function $B \mapsto \lambda(B)$
in the space of continuous maps $B: M\to\SL(2)$ equipped with the $C^0$ topology.
Moreover, $\lambda(A)>0$ if and only if there exists some periodic point $p\in M$ such that
$A^{\per(p)}(p)$ is a hyperbolic matrix.
In that case, $\lambda(\cdot)$ is bounded from zero for all continuous cocycles on a
$C^0$-neighborhood of $A$.
\end{main}

This applies, in particular, in the setting of Young~\cite{You93} and thus Theorem~\ref{t.main2}
contains a much stronger version of the main result in there:
the Lyapunov exponent is $C^0$-continuous at every $C^2$ cocycle in the isotopy class;
moreover, it is non-zero if and only if the cocycle is hyperbolic on some periodic orbit
(an open and dense condition).

All the cocycles we consider are of the form $\hF_B(\hx.v)=(\hf(\hx),B\circ\pi(\hx)v)$ for
some continuous $B:M\to\SL(2)$ and so they all have (trivial) $s$-holonomy.
Thus the notion of $s$-invariant measure, as defined in \eqref{eq.s-invariant}, makes
sense for such cocycles. In Section~\ref{s.invariant} we study certain properties of such
measures. We do not assume the cocycle to be $u$-bunched, and so the conclusions apply for
both theorems. In Section~\ref{s.main1} we deduce Theorem~\ref{t.main1}.

The proof of Theorem~\ref{t.main2} is similar, but we have to deal with the fact that
$u$-holonomies need not exist, since we do not assume $u$-bunching.
The first step, in Section~\ref{s.nobunching}, is to explain what we mean by a
$u$-invariant measure and a $u$-state. Next, we need a suitable version of the
Invariance Principle of \cite{Extremal,BGV03}. This we prove in Section~\ref{s.invariance},
using ideas from \cite{TY16}. In Section~\ref{s.states}, we check that the assumptions
ensure that there are no $su$-states. In Section~\ref{s.main2} we wrap up the proof.

\section{$s$-invariant measures}\label{s.invariant}

Let $\cM^s$ be the space of probability measures on $\hM\times\pR^2$ that project down to
$\hmu$ and are $s$-invariant.
Let $\cM$ be the space of probability measures on $M\times\pR^2$ that project down to $\mu$.
Both spaces are equipped with the weak$^*$ topology.
Consider the map $\Psi:\cM\to\cM^s$ defined as follows:
given any $m\in\cM$ and a disintegration $\{m_x:x\in M\}$ along the fibers $\{x\}\times\pR^2$,
let $\hm=\Psi(m)$ be the measure on $\hM\times\pR^2$ that projects down to $\hmu$ and
whose conditional probabilities $\hm_\hx$ along the fibers $\{\hx\} \times\pR^2$
are given by
$$
\hm_\hx=m_{\pi(\hx)}.
$$

It is clear from the definition that $\hm\in\cM^s$.
Moreover, if $\{m'_x:x\in M\}$ is another disintegration of $m$ then, by essential uniqueness of
the disintegration, $m_x=m'_x$ for $\mu$-almost every $x$. Recalling that $\hmu$ is the lift of
$\mu$, this implies that $m_{\pi(\hx)} = m'_{\pi(\hx)}$ for $\hmu$-almost every $\hx$.
Thus $\hm$ does not depend on the choice of the disintegration.
This shows that $\Psi$ is well-defined. We are going to prove:

\begin{proposition}\label{p.homeomorphism}
$\Psi:\cM\to\cM^s$ is a homeomorphism.
\end{proposition}

\begin{proof}
We use the fact that $\hmu$ has local product structure (see \cite[Section~2.2]{BoV04}).
For each $\hp\in\hM$, let $p=\pi(\hp)$ and consider the neighborhood $\hV_\hp = \pi^{-1}(B(p,\rho))$.
We may identify $\hV_\hp$ to the product
$$
B(p,\rho) \times \pi^{-1}(p) = W^u_\loc(\hp) \times W^s_\loc(\hp)
$$
through a homeomorphism, in such a way that $\pi$ becomes the projection to the first coordinate.
Local product structure gives that the restriction of $\hmu$ may be written as
$$
\hmu\mid\hV_\hp = \rho \, \hmu^u \times \hmu^s,
$$
where $\rho:\hV_\hp\to (0,\infty)$ is a continuous function, $\hmu = \mu \mid B(p,\rho)$ and $\mu^s$ is a probability measure on
$W_\loc^s(\hp)$.

\begin{lemma}\label{l.decomposition}
For any $m\in\cM$, the measure $\hm=\Psi(m)$ satisfies
$$
\hm \mid \hV_\hp \times \pR^2 = \rho \, (m \mid B(p,\rho)) \times \hmu^s
\text{ for any } \hp \in\hM.
$$
\end{lemma}

\begin{proof}
Given any bounded measurable function $g:\hV_\hp\times\pR^2\to\RR$,
$$
\begin{aligned}
\int_{\hV_\hp \times \pR^2} g \, d\hm
& = \int_{\hV_\hp} \int_{\pR^2} g(\hx,v) \, d\hm_{\hx}(v) \, d\hmu(\hx)\\
& = \int_{W_\loc^s(\hp)} \int_{W_\loc^u(\hp)} \int_{\pR^2} g(x,\xi,v) \, d\hm_{(x,\xi)}(v) \, \rho(x,\xi) \, d\hmu^u(x) \, d\hmu^s(\xi).
\end{aligned}
$$
Since $\hm_{(x,\xi)} = m_x$ for every $x\in M$, by definition, it follows that
$$
\begin{aligned}
\int_{\hV_\hp \times \pR^2} g \, d\hm
& = \int_{W_\loc^s(\hp)} \int_{W_\loc^u(\hp)} \int_{\pR^2}
g(x,\xi,v) \rho(x,\xi) \, dm_x (v) \, d\hmu^u(x) \, d\hmu^s(\xi)\\
& = \int_{W_\loc^s(\hp)} \int_{W_\loc^u(\hp) \times \pR^2}
g(x,\xi,v) \rho(x,\xi) \, dm(x,v) \, d\hmu^s(\xi).
\end{aligned}
$$
This proves the claim.
\end{proof}

Let us prove that $\Psi$ is continuous, that is, that given any sequence $(m_n)_n$
converging to some $m$ in $\cM$, we have
\begin{equation}\label{eq.donverges}
\int_{\hM\times\pR^2} g \, d\Psi(m_n) \to \int_{\hM\times\pR^2} g \, d\Psi(m)
\end{equation}
for every continuous  function $g:\hM\times\pR^2\to\RR$. It is no restriction to suppose that the
support of $g$ is contained in $\hV_\hp$ for some $\hp\in\hM$, for every continuous function
is a finite sum of such functions.
Then, by Lemma~\ref{l.decomposition},
$$
\int_{\hM\times\pR^2} g \, d\Psi(m_n)
=  \int_{W_\loc^u(\hp) \times \pR^2} \int_{W_\loc^s(\hp)}
g(x,\xi,v) \rho(x,\xi)  \, d\hmu^s(\xi) \, dm_n(x,v).
$$
Our hypotheses ensure that $G(x,v) = \int_{W_\loc^s(\hp)}  g(x,\xi,v) \rho(x,\xi)  \, d\hmu^s(\xi)$
defines a continuous function. Hence, the assumption that $(m_n) \to m$ implies that
$$
\begin{aligned}
\int_{\hM\times\pR^2} g \, d\Psi(m_n)
 & = \int_{W_\loc^u(\hp) \times \pR^2} G(x,v) \, dm_n(x,v) \to \\
& \qquad \to \int_{W_\loc^u(\hp) \times \pR^2} G(x,v) \, dm(x,v)
= \int_{\hM\times\pR^2} g \, d\Psi(m),
\end{aligned}
$$
as claimed. We are left to proving that $\Psi$ is a bijection.

Surjectivity is clear: given $\hm\in\cM^s$ take $m$ to be the probability measure on
$M\times\pR^2$ that projects down to $\mu$ and whose conditional probabilities along the vertical
fibers $\{x\}\times\pR^2$ are given by $m_x = \hm_\hx$ for any $\hx\in\pi^{-1}(x)$.
This is well defined, by \eqref{eq.s-invariant}, and it is clear from the definition
that $\Psi(m) = \hm$.
Injectivity is a consequence of Lemma~\ref{l.decomposition}.
Indeed, if $\Psi(m_1)=\Psi(m_2)$ then
$$
\int_{X \times V} \int_{W^s_\loc(\hp)} \rho(x,\xi) \, d\hmu^s(\xi) \, dm_1(x,v)
=
\int_{X \times V} \int_{W^s_\loc(\hp)} \rho(x,\xi) \, d\hmu^s(\xi) \, dm_2(x,v)
$$
for any $\hp\in\hM$ and any $X\times V\subset B(p,\rho) \times\pR^2$. This implies that
$$
Hm_1\mid B(p,\rho)=Hm_2\mid B(p,\rho),
\quad\text{where $H(x) = \int_{W^s_\loc(\hp)} \rho(x,\xi) \, d\hmu^s(\xi)$.}
$$
Noting that $H$ is positive, it follows that the restrictions of $m_1$ and $m_2$ to
$B(p,\rho)$ coincide, for every $p\in M$. Thus $m_1=m_2$.
\end{proof}

\begin{corollary}\label{c.compact}
$\cM^s$ is non-empty, convex and compact.
\end{corollary}

\begin{proof}
Convexity is obvious and the other claims follow directly from Proposition~\ref{p.homeomorphism},
since $\cM$ is non-empty and compact.
\end{proof}

Let $(B_n)_n$ be a sequence of maps converging uniformly to some $B$ in the space of continuous maps
$\cM\to\SL(2)$, and $(\hm_n)_n$ be a sequence of probability measures on $\hM$ converging in the
weak$^*$ topology to some probability measure $m$.

\begin{corollary}\label{c.closed}
If $\hm_n$ is an $s$-state of $B_n$ for every $n$ then $\hm$ is an $s$-state of $B$.
\end{corollary}

\begin{proof}
It follows from Corollary~\ref{c.compact} that $\hm\in\hM^s$. It is clear that $\hm$ is
$\hF_B$-invariant, because $m_n$ is $\hF_{B_n}$-invariant for every $n$ and $\hF_{B_n}$
converges uniformly to $\hF_B$.
\end{proof}

\section{Proof of Theorem~\ref{t.main1}}\label{s.main1}

If $\lambda(A)=0$ then, by the Invariance Principle (\cite[Theorem~D]{Extremal}, \cite[Th\'eor\`eme~1]{BGV03}),
every $\hF_A$-invariant probability measure $\hm$ that projects down to $\hmu$ is an $su$-state.
Thus, the hypothesis that $A$ has no $su$-states implies that $\lambda(A)>0$.

We are left to proving that $A$ is a continuity point for the Lyapunov exponent.
Define (here $v$ denotes both a direction and any non-zero vector in that direction)
$$
\phi_B:\hM\times\pR^2 \to \RR,
\quad
\phi_B(\hx,v) = \log \frac{\|\hB(\hx)v\|}{\|v\|}.
$$

\begin{proposition}\label{p.realization}
Every $B:M\to\SL(2)$ in a $C^0$-neighborhood of $A$ admits some $s$-state $\hm_B$
such that
$$
- \lambda(B) = \int_{\hM\times\pR^2} \phi_B \, d\hm.
$$
\end{proposition}

\begin{proof}
First, suppose that $\lambda(B)=0$. For every $(\hx,v)\in\hM\times\pR^2$ and $n\ge 1$,
$$
\sum_{j=0}^{n-1} \phi_B(\hF_B^j(\hx,v)) = \log \frac{\|\hB^n(\hx)v\|}{\|v\|}
\in \Big[-\log\|\hB^n(\hx)^{-1}\|, \log \|\hB^n(\hx)\| \Big].
$$
We also have that, or $\hmu$-almost every $\hx\in M$,
$$
\lim_n \frac 1n \log \|\hB^n(\hx)\|
= \lim_n \frac 1n \log \|\hB^n(\hx)^{-1}\|
= \lambda(B).
$$
Thus, for any $\hF_B$-invariant measure $\hm$ that projects down to $\hmu$,
$$
\lim_n \frac 1n \sum_{j=0}^{n-1} \phi_B(\hF_B^j(\hx,v)) = 0
\text{ for $\hm$-almost every $(\hx,v)$.}
$$
By the ergodic theorem, this implies that
$$
\int_{\hM\times\pR^2} \phi_B \, d\hm = 0 = \lambda(B)
$$
and so every $\hF_B$-invariant measure $\hm$ that projects down to $\hmu$ satisfies the
conclusion of the lemma.

Now suppose that $\lambda(B)>0$. By the theorem of Oseledets~\cite{Ose68}, there exists
an $\hF_B$-invariant splitting $\RR^2 = E^u_\hx \oplus E^s_\hx$ defined at $\hmu$-almost
every point $\hx$ and such that
\begin{equation}\label{eq.splitting}
\begin{aligned}
& \lim_{n\to\pm\infty} \frac 1n \log\|\hB^n(\hx)v\| =  \lambda(B) \text{ for non-zero } v \in E^u_\hx
\text{ and }\\
& \lim_{n\to\pm\infty} \frac 1n \log\|\hB^n(\hx)v\| = - \lambda(B)\text{ for non-zero } v \in E^s_\hx.
\end{aligned}
\end{equation}
Let $\hm$ be the probability measure on $\hM\times\pR^2$ that projects down to $\hmu$ and
whose conditional probabilities along the fibers $\{\hx\} \times\pR^2$ are given by the Dirac
masses at $E^s_\hx$. Then $\hm$ is an $s$-state: it is clear that it is $\hF_B$-invariant;
to check that it is $s$-invariant, just note that the subspace $E^s_\hx$ depends only on the
forward iterates, and so it is constant on each $\pi^{-1}(x)$.
Moreover, by the ergodic theorem and \eqref{eq.splitting},
$$
\begin{aligned}
\int_{\hM\times\pR^2} \phi_B \, d\hm
& = \int_{\hM\times\pR^2} \lim_n \frac 1n \sum_{j=0}^{n-1} \phi_B \circ \hF_B^j \, d\hm \\
& = \int_{\hM} \int_{\pR^2}\lim_n \frac 1n \log  \frac{\|\hB^n(\hx)v\|}{\|v\|} \, d\delta_{E^s_\hx}(v) \, d\hmu(\hx)\\
& = \int_{\hM} - \lambda(B) \, d\hmu(\hx) = - \lambda(B).
\end{aligned}
$$
This completes the proof.
\end{proof}

\begin{lemma}\label{l.sus}
If $A$ has no $su$-states then it has exactly one $s$-state.
\end{lemma}

\begin{proof}
Existence is contained in Proposition~\ref{p.realization}. To prove uniqueness, we argue as follows.
Let $\hm$ be any $s$-state. As observed before, the hypothesis implies that $\lambda(A)>0$.
Let $\RR^2 = E^u_\hx \oplus E^s_\hx$ be the Oseledets invariant splitting, defined at $\hmu$-almost
every point $\hx$. Let $\hm^u$ and $\hm^s$ be the probability measures on $\hM\times\pR^2$ that
project down to $\hmu$ and whose conditional probabilities along the fibers $\{\hx\} \times\pR^2$
are the Dirac masses at $E^u_\hx$ and $E^s_\hx$, respectively.
Then $\hm^u$ is a $u$-state, $\hm^s$ is an $s$-state and every $\hF_A$-invariant probability
measure is a convex combination of $\hm^u$ and $\hm^s$ (compare~\cite[Lemma~6.1]{Extremal}).
Then, keeping in mind that $\hmu$ is ergodic, there is $\alpha\in [0,1]$ such that
$\hm = \alpha \hmu^u + (1-\alpha) \hmu^s$. If $\alpha>0$, we may write
$$
\hm^u = \frac{1}{\alpha} \hm + \big(1-\frac{1}{\alpha}) \hm^s
$$
and, as $\hm$ and $\hm^s$ are $s$-states, it follows that $\hm^u$ is an $s$-state.
Since $\hm^u$ is also a $u$-state, this contradicts the hypothesis.
Thus $\alpha=0$, that is, $m=m^s$.
\end{proof}

Theorem~\ref{t.main1} is an easy consequence. Indeed, we already know that $\lambda(A)>0$.
Consider any sequence $A_k: M\to\SL(2)$, $k\in\NN$ converging to $A$ in the $C^0$ topology.
By Proposition~\ref{p.realization}, for each $k$ one can find some $s$-state $\hm_k$ for $A_k$
such that
$$
- \lambda(A_k) = \int_{\hM\times\pR^2} \phi_{A_k} \, d\hm_k.
$$
Up to restricting to a subsequence, we may suppose that $(\hm_k)_k$ converges to some
probability measure $\hm$ in the weak$^*$ topology.
By Corollary~\ref{c.closed}, $\hm$ is an $s$-state for $A$. Moreover, since $\phi_{A_k}$
converges uniformly to $\phi_{A}$,
\begin{equation}\label{eq.lambda}
\lim_k - \lambda(A_k) = \int_{\hM\times\pR^2} \phi_{A} \, d\hm.
\end{equation}
By Proposition~\ref{p.realization} and Lemma~\ref{l.sus}, the right-hand side is equal to $\lambda(A)$.
This proves continuity of the Lyapunov exponent in the $C^0$ topology.

\begin{remark}\label{r.sus}
The converse to Lemma~\ref{l.sus} is true when $\lambda(A)>0$.
\end{remark}

\section{$u$-states without $u$-bunching}\label{s.nobunching}

Next we prove Theorem~\ref{t.main2}.
Initially, suppose that $0 \le \lambda(A) < \log\sigma$. Then the cocycle is ``nonuniformly $u$-bunched,''
in a sense that was exploited before, in \cite[Sections~2.1 and~2.2]{Almost}. Using those methods,
one gets that (compare \cite[Proposition~2.5]{Almost})
$$
h_{\hx,\hy}^u = \lim_n \hA^n\big(\hf^{-n}(\hy)\big) \hA^n\big(\hf^{-n}(\hx)\big)^{-1}
$$
exists for $\hmu$-almost every $\hx$ and any $\hy\in W^u_\loc(\hx)$.
Then we define a probability measure $\hm$ on $\hM\times\pR^2$ to be \emph{$u$-invariant}
if it admits a disintegration $\{\hm_\hx: \hx \in\hM\}$ along the fibers $\{\hx\}\times\pR^2$ such that
\begin{equation}\label{eq.u-invariant2}
\big(h_{\hx,\hy}^u\big)_*\hm_{\hx} = \hm_{\hy}
\text{ for $\hmu$-almost every $\hx$ and any $\hy\in W^u_\loc(\hx)$.}
\end{equation}
As before, a \emph{$u$-state} is an $\hF_A$-invariant probability measure which is $u$-invariant.

When $\lambda(A) \ge \log \sigma$ we have to restrict ourselves to the subclass of $\hF_A$-invariant
probability measures whose center Lyapunov exponent is strictly less than $\sigma$.
More precisely, we consider only $\hF_A$-invariant probability measures $\hm$ such that
\begin{equation}\label{eq.boundexp}
\lim_n \frac 1n \log \|D\hA^n(\hx)v\| \le c < \log\sigma
\text{ for $\hm$-almost every $(\hx,v)\in\hM\times\pR^2$,}
\end{equation}
where $D\hA(\hx)v$ denotes the derivative of the projective map $\hA(\hx):\pR^2\to\pR^2$.

\begin{remark}\label{r.derivative}
The following elementary bound will be useful:
$$
\|\hA(\hx)\|^{-1}\|\hA(\hx)^{-1}\|^{-1}
\le \frac{\|D\hA(\hx)v\|}{\|v\|}
\le \|\hA(\hx)\|\|\hA(\hx)^{-1}\|
\text{ for every $\hx$.}
$$
\end{remark}

In the next result we use the fact that the natural extension of $f$ admits a $C^\H$ realization:
there exist a $C^\H$ embedding $g:U \to U$ defined on some open subset $U$ of an Euclidean space,
and a topological embedding $\iota:\hM \to U$ with $g(\iota(\hM)) = \iota(\hM)$ and
$g \circ\iota = \iota \circ \hf$. A proof is given in Appendix~\ref{a.natural}.
Identifying $\hM$ with $\iota(\hM)$ we may view $\hf$ as a restriction of $g$,
and so we may apply Pesin theory to it.

\begin{proposition}\label{p.Pesin1}
If $\hm$ satisfies \eqref{eq.boundexp} then for $(\hx,v)$ in a full $\hm$-measure
subset $\Lambda$  of $\hM\times\pR^2$ there exists a $C^1$ function
$\psi_{\hx,v}:W^u_\loc(\hx)\to\pR^2$ depending measurably on $(\hx,v)$ such that
$\psi_{\hx,v}(\hx)=v$ and the graphs
$\cW^u_\loc(\hx,v)=\{(\hy,\psi_{\hx,v}(\hy)):\hy\in W^u_\loc(\hx)\}$ satisfy
\begin{enumerate}
\item[(a)] $\hF^{-1}(\cW^u_\loc(\hx,v)) \subset\cW^u_\loc(\hF^{-1}(\hx,v))$ for every $(\hx,v)\in\Lambda$;
\item[(b)] $d(\hF^{-n}(\hx,v),\hF^{-n}(\hy,w)) \to 0$ exponentially fast for any $(\hy,w)\in \cW^u_\loc(\hx,v)$.
\end{enumerate}
\end{proposition}

\begin{proof}
The assumption ensures that there exists $\hm$-almost everywhere an Oseledets strong-unstable
subspace $\hE^u_{\hx,v}\subset T_{\hx}U \times \RR^2$ that is a graph over the unstable
direction $E_\hx\subset T_\hx U$ of $g$. Then, by Pesin theory, there exists $\hm$-almost
everywhere a $C^1$ embedded disk $\tW^u(\hx,v)$ tangent to  $\hE^u_{\hx,v}$ and such that
$$
\hF^{-n}(\hy,w) \in \tW^u_\loc(\hF^{-n}(\hx,v))
\quand
d(\hF^{-n}(\hx,v),\hF^{-n}(\hy,w)) \le \sigma^{-n}
$$
for every $n\ge 0$ and $(\hy,w)\in \tW^u_\loc(\hx,v)$. This also implies that $\tW_\loc^u(\hx,v)$
is a $C^1$ graph over a neighborhood of $\hx$ inside $W^u(\hx)$. While the radius $r(\hx)$ of this
neighborhood need not be bounded from zero, in principle, Pesin theory also gives that it decreases
sub-exponentially along orbits:
$$
\lim_n \frac 1n \log r(\hf^{-n}(\hx)) = 0.
$$
On the other hand, the size of $\hf^{-n}(W^u_\loc(\hx))$ decreases exponentially fast (faster than
$\sigma^{-n}$). Thus, the projection of $\tW^u(\hF^{-n}(\hx,v))$ contains $\hf^{-n}(W^u_\loc(\hx))$
for any large $n$. Then $\hF^n(\tW^u(\hF^{-n}(\hx,v)))$ is a $C^1$ graph whose projection contains
$W^u_\loc(\hx)$. Take $\cW^u_\loc(\hx,v)$ to be the part of that graph that lies over $W^u_\loc(\hx)$.
It is clear from the construction that conditions (a) and (b) in the statement are satisfied.
\end{proof}

Denote $\Lambda_\hx=\Lambda \cap (\{\hx\}\times\pR^2)$ for each $\hx\in\hM$.
We say that $\hm$ is \emph{$u$-invariant} if it admits a disintegration $\{\hm_\hx:\hx\in\hM\}$
along the fibers $\{\hx\}\times\pR^2$ such that
\begin{equation}\label{eq.u-invariant3}
\big(h_{\hx,\hy}^u\big)_*\hm_{\hx} = \hm_{\hy}
\text{ for $\hmu$-almost every $\hx$ and any $\hy\in W^u_\loc(\hx)$,}
\end{equation}
where $h^u_{\hx,\hy}:\Lambda_\hx \to \{\hy\}\times\pR^2$ is defined by
$h^u_{\hx,\hy}(\hx,v) = (\hy,\psi_{(\hx,v)}(\hy))$. Note that $\hm_\hx(\Lambda_\hx)=1$ for
$\hmu$-almost every $\hx$, because $\hm(\Lambda)=1$.
By definition, a \emph{$u$-state} is an $\hF_A$-invariant probability measure $\hm$ that
satisfies \eqref{eq.boundexp} and is $u$-invariant.

\section{A new $u$-invariance principle}\label{s.invariance}

Here we prove the following form of the Invariance Principle, where the main novelty is that no
$u$-bunching is assumed:

\begin{theorem}\label{t.principle}
Every $\hF_A$-invariant probability measure $\hm$ satisfying
\begin{equation}\label{eq.negatexp}
\lim_n \frac 1n \log \|D\hA^n(\hx)v\| \le 0
\text{ for $\hm$-almost every $(\hx,v)\in\hM\times\pR^2$,}
\end{equation}
is a $u$-state.
\end{theorem}

We are going to extend to our setting an approach introduced by Tahzibi, Yang~\cite{TY16} for
bunched cocycles. This is based on the notion of partial entropy, which may be defined
as follows (see \cite{Led84a,Yang-partial} for more information).

Let $\cR$ be a Markov partition of $\hf$ with diameter small enough that $\cR(\hx) \subset \hV_\hx$
for every $\hx \in \hM$, where $\cR(\hx)$ denotes the element of $\cR$ that contains $\hx$.
(Actually, elements of $\cR$ may intersect along their boundaries but, since the boundaries
are nowhere dense and have zero $\hmu$-measure, we may ignore the trajectories through them.)
Let $\xi^u(\hx)\subset\hM$ be the connected component of $\cR(\hx) \cap W_\loc^u(\hx)$ that contains $\hx$.
For $v\in\pR^2$ such that $(\hx,v)\in\Lambda$, let $\Xi^u(\hx,v)$ be the connected component of
$(\cR(x) \times\pR^2) \cap \cW^u_\loc(\hx,v)$ that contains $(\hx,v)$.

The family $\xi^u$ is an \emph{adapted partition} for $(\hf,\hmu)$: its elements are pairwise
disjoint and, for $\hmu$-almost every $\hx$,
\begin{itemize}
\item $\hf^{-1}(\xi^u(\hx)) \subset \xi^u(\hf^{-1}(\hx))$ and
\item $\xi^u(\hx)$ contains a neighborhood of $\hx$ inside $W^u(\hx)$.
\end{itemize}
Analogously, $\Xi^u$ is an adapted partition for $(\hF,\hm)$.
The corresponding \emph{partial entropies} are defined by
\begin{equation}\label{eq.entropies}
h_\hmu(\hf,W^u)=H_{\hmu}(\hf^{-1}\xi^u \mid \xi^u)
\quand
h_\hm(\hF_A,\cW^u)=H_{\hm}(\hF_A^{-1}\Xi^u \mid \Xi^u).
\end{equation}

\subsection{$c$-invariant measures}\label{s.cinvariance}

Let $\{\hmu^u_\hx:\hx \in \hM\}$ and $\{\hm^u_{\hx,v}:(\hx,v) \in \hM\times\pR^2\}$ be
disintegrations of, respectively, $\hmu$ relative to $\xi^u$ and $\hm$ relative to $\Xi^u$.
Let $p:\hM \times \pR^2\to \hM$ be the canonical projection. We call $\hm$ \emph{$c$-invariant} if
\begin{equation}\label{eq.cinvariant1}
(h^c_{\hx,v,w})_* \hm^u_{\hx,v}=\hm^u_{\hx,w}
\text{ for $\hm$-almost  every $(\hx,v)$ and $(\hx,w)$,}
\end{equation}
where $h^c_{\hx,v,w}:\Xi^u(\hx,v)\to\Xi^u(\hx,w)$ is the bijection defined by $p\circ h^c_{\hx,v,w}=p$.
Equivalently, the measure $\hm$ is $c$-invariant if
\begin{equation}\label{eq.cinvariant2}
p_*(\hm^u_{\hx,v})=\hmu^u_{\hx}
\text{ for $\hm$-almost  every $(\hx,v)$.}
\end{equation}

\begin{proposition}\label{p.cuequivalent}
The measure $\hm$ is $u$-invariant if and only if it is $c$-invariant.
\end{proposition}

\begin{proof}
Let us start with a model: let $\nu$ be a probability measure on a product $X \times Y$
of two separable metric spaces, and let $\{\nu^1_y:y\in Y\}$ and $\{\nu^2_x:x\in X\}$ be
disintegrations of $\nu$ relative to the partition into horizontals $X\times\{y\}$ and
the partition into verticals $\{x\}\times Y$, respectively.
We call $\nu$ \emph{$v$-invariant} (respectively, \emph{$h$-invariant}) if the disintegrations
may be chosen such that $\nu_y^1$ is independent of $y$
(respectively, $\nu^2_x$ is independent of $x$).

\begin{lemma}\label{l.hv}
$\nu$ is $v$-invariant if and only if it is $h$-invariant.
\end{lemma}

\begin{proof}
Suppose that $\nu$ is $v$-invariant and let $\nu^1$ be such that $\nu^1_y=\nu^1$ for every $y$.
Let $\nu^2$ be the quotient of $\nu$ relative to the horizontal partition or, equivalently,
the projection of $\nu$ to the second coordinate. Then, by the definition of disintegration,
$$
\nu = \nu^1 \times \nu^2.
$$
This implies that $\nu^1$ is the projection of $\nu$ to the first coordinate and $\nu_x^2=\nu_2$
defines a disintegration of $\nu$ relative to the vertical partition.
In particular, $\nu$ is $h$-invariant. The proof that $h$-invariance implies $v$-invariance is identical.
\end{proof}

To deduce the proposition we only have to reduce its setting to that of Lemma~\ref{l.hv}.
Consider the partitions $\Xi^c$ and $\Xi^{cu}$ of $\hM\times\pR^2$ defined by
$$
\Xi^c(\hx,v) = p^{-1}(\hx)
\quand
\Xi^{cu}(\hx,v) = p^{-1}(\xi^u(\hx)).
$$
Let $\{\hm^c_{\hx,v}: (\hx,v) \in \hM\times\pR^2\}$ and $\{\hm^{cu}_{\hx,v}: (\hx,v) \in \hM\times\pR^2\}$
be disintegrations of $\hm$ relative to $\Xi^c$ and $\Xi^{cu}$, respectively.
Both $\Xi^c$ and $\Xi^u$ refine $\Xi^{cu}$. So, by essential uniqueness of the disintegration,
\begin{enumerate}
\item[(i)]
$\{\hm^u_{\hy,w}: (\hy,w)\in\Xi^{cu}(\hx,v)\}$ is a disintegration of $\hm^{cu}_{\hx,v}$ with respect
to the partition $\Xi^u \mid \Xi^{cu}(\hx,v)$ and
\item[(ii)]
$\{\hm^c_{\hy,w}: (\hy,w)\in\Xi^{cu}(\hx,v)\}$ is a disintegration of $\hm^{cu}_{\hx,v}$ with respect
to the partition $\Xi^c\mid \Xi^{cu}(\hx,v)$,
\end{enumerate}
for $\hm$-almost every $(\hx,v)$. This will be used a few times in the following.

Now consider the map
$$
\Psi_{\hx,v}: \Xi^{cu}(\hx,v) \to \xi^u(\hx) \times \pR^2,
\quad \Phi_\hx(\hy,w) = (\hy,z)
$$
where $z$ is such that $(\hx,z)$ is the point where $\Xi^u(\hy,w)$ intersects $\Xi^c(\hx,v)$.
Since $\Lambda$ has full $\hm$-measure, $\Psi_{\hx,v}$ is defined $\hm^{cu}_{\hx,v}$-almost everywhere
for $\hm$-almost every $(\hx,v)$. Clearly, it is an invertible measurable map sending
atoms of $\Xi^u\mid \Xi^{cu}(\hx,v)$ to horizontals $\xi^u(\hx) \times \{z\}$ and
atoms of $\Xi^c \mid \Xi^{cu}(\hx,v)$ to verticals $\{\hy\} \times \pR^2$.

Identifying $\Xi^{cu}(\hx,v)$ to $\xi^u(\hx) \times \pR^2$ through $\Psi_{\hx,v}$, (i) and (ii) above
correspond to disintegrations of $\hm_{\hx,v}$ relative to the horizontal partition and the vertical
partition, respectively.
Moreover, $s$-invariance and $u$-invariance translate to $v$-invariance and $h$-invariance, respectively.
Thus the claim follows from Lemma~\ref{l.hv}.
\end{proof}

\subsection{A criterion for $c$-invariance}\label{s.criterion}

Note that $h_\hmu(\hf) \le h_\hm(\hF_A)$, because $(\hf,\hmu)$ is a factor of $(\hF_A, \hm)$.
For the partial entropies the inequality goes the opposite way:

\begin{proposition}\label{p.criterion}
$h_\hm(\hF_A,\cW^u) \leq h_\hmu(\hf, W^u)$ and the equality holds if and only if $\hm$ is $c$-invariant.
\end{proposition}

\begin{proof}
Keep in mind that $\xi^u \prec \hf^{-1}\xi^u$ and $\Xi^u \prec \hF_A^{-1}\Xi^u$. By definition,
\begin{equation}\label{eq.h1}
\begin{aligned}
h_\hmu(\hf,W^u)
= H_\hmu (\hf^{-1}\xi^u \mid \xi^u)
& = \int H_{\hmu^u_\hx}(\hf^{-1}\xi^u) \, d\hmu(\hx) \text{ where}\\
H_{\hmu^u_\hx} (\hf^{-1}\xi^u)
& = \int - \log \hmu^u_\hx(\hf^{-1}\xi^u(\hy)) \, d\hmu^u_\hx(\hy),
\end{aligned}
\end{equation}
and similarly for $h_\hm(\hF_A,\cW^u)$ and $\Xi^{u}$.

\begin{lemma}\label{l.local}
For $\hm$-almost every $(\hx,v)\in\hM\times\pR^2$,
\begin{enumerate}
\item[(a)] $H_{\hm^{cu}_{\hx,v}}(\hF_A^{-1} \Xi^{u} \mid \Xi^u) \leq H_{\hmu^u_\hx}(\hf^{-1}\xi^u)$ and
\item[(b)] the equality holds if and only if $\hm^u_{\hx,v}(\hat{F}_A^{-1}\Xi^u(\hy,w))=\hmu^u_\hx(\hf^{-1}\xi^u(\hy))$
for $\hm^{cu}_{\hx,v}$-almost every $(\hy,w)\in\Xi^{cu}(\hx,v)$.
\end{enumerate}
\end{lemma}

\begin{proof}
Since $\hmu=p_*\hm$ and $\Xi^{cu}(\hx,v)=p^{-1}(\xi^u(\hx))$, essential uniqueness of disintegrations gives that
$\hmu^u_\hx = p_*(\hm^{cu}_{\hx,v})$ for $\hm$-almost every $(\hx,v)$. Thus,
$$
\begin{aligned}
H_{\hmu^u_\hx}(\hf^{-1}\xi^u)
& = \int - \log \hmu^u_\hx(\hf^{-1}\xi^u(\hy)) \, d\hmu^{u}_\hx(\hy) \\
& = \int - \log \hm^{cu}_{\hx,v}(\hF_A^{-1}\Xi^{cu}(\hy,w)) \, d\hm^{cu}_{\hx,v}(\hy,w)
= H_{\hm^{cu}_{\hx,v}}(\hF_A^{-1}\Xi^{cu})
\end{aligned}
$$
for $\hm$-almost every $(\hx,v)$. Moreover, using the relation
$\hF^{-1}_A \Xi^{cu} \vee \Xi^u = \hF^{-1}_A \Xi^u$,
$$
H_{\hm^{cu}_{\hx,v}}(\hF^{-1}_A \Xi^{cu})
\geq H_{\hm^{cu}_{\hx,v}}(\hF^{-1}_A \Xi^{cu} \mid \Xi^{u})
= H_{\hm^{cu}_{\hx,v}}(\hF^{-1}_A \Xi^u \mid \Xi^{u}).
$$
This proves claim (a). Moreover, the equality holds if and only if the partitions
$\hF^{-1}_A \Xi^{cu}$ and $\Xi^{u}$ are independent relative to $\hm^{cu}_{\hx,v}$,
that is,
$$
\hm^u_{\hx,v}(\hF^{-1}_A \Xi^{cu}(\hy,w))
= \hm^{cu}_{\hx,v}(\hF^{-1}_A \Xi^{cu}(\hy,w))
\text{ for $\hm^{cu}_{\hx,v}$-almost every $(\hy,w)$.}
$$
By the previous observations, this is equivalent to
$$
\hm^u_{\hx,v}(\hF^{-1}_A \Xi^u(\hy,w))
= \hmu^u_\hx(\hf^{-1} \xi^u(\hy))
\text{ for $\hm^{cu}_{\hx,v}$-almost every $(\hy,w)$,}
$$
as claimed in (b).
\end{proof}

Similarly to \eqref{eq.h1}, we have
$H_{\hm^{cu}_{\hx,v}}(\hF^{-1}_A\Xi^u \mid \Xi^u) = \int H_{\hm^u_{\hy,w}}(\hF_A^{-1}\Xi^u) \, d\hm^{cu}_{\hx,v}(\hy,v)$.
So, integrating the inequality in part (a) of the lemma,
$$
\begin{aligned}
H_\hm(\hF_A^{-1} \Xi^u \mid \Xi^u)
& = \int H_{\hm^u_{\hx,v}}(\hF_A^{-1} \Xi^u) \, d\hm(\hx,v)
= \int H_{\hm^{cu}_{\hx,v}}(\hF_A^{-1} \Xi^u \mid \Xi^u) \, d\hm(\hx,v) \\
& \le \int H_{\mu^u_\hx}(\hf^{-1}\xi^u) \, d\hmu(\hx)
= H_\hmu(\hf^{-1}\xi \mid \xi^u).
\end{aligned}
$$
Moreover, the equality holds if and only if
$\hm^u_{\hx,v}(\hF_A^{-1} \Xi^u(\hx,v))=\hmu^u_\hx(\hf^{-1}\xi^u(\hx))$ for $\hm$-almost every $(\hx,v)$.
In other words, the equality holds if and only if $p_*\hm^u_{\hx,v}=\hmu^u_\hx$ restricted to the
$\sigma$-algebra generated by $\hF_A^{-1}\Xi^u$.

Replacing $\hF_A$ by any iterate $\hF^n_A$, and noting that
$$
h_\hm(\hF_A^n,\cW^u) = n h_\hm(\hF_A,\cW^u) \quand  h_\hmu(\hf^n,W^u) = n h_\hmu(\hf,W^u),
$$
we get that the equality holds if and only if $p_*\hm^u_{\hx,v}=\hmu^u_\hx$ restricted to the $\sigma$-algebra
generated by $\hF^{-n}_A \Xi^u$. Since $\cup_n \hF^{-n}_A \Xi^u$ generates the Borel $\sigma$-algebra of every
$\Xi^u(\hx,v)$, this is the same as saying that $p_*\hm^u_{\hx,v}=\hmu^u_\hx$ for $\hm$-almost every $(\hx,v)$,
that is, that $\hm$ is $c$-invariant.
\end{proof}

\subsection{Proof of Theorem~\ref{t.principle}}\label{ss.equalentropy}

The hypothesis \eqref{eq.negatexp} ensures that the Lyapunov exponents of $\hm$ along the center
fibers $\{\hx\}\times\pR^2$ are non-positive. Then
$$
h_\hm(\hF_A) = h_\hm(\hF_A , \cW^u)
$$
(see \cite[Corollary 5.3]{LeY85A}). Similarly, $h_\hmu(\hf) = h_\hmu(\hf,W^u)$.
Moreover, $h_\hm(\hF_A) \geq h_\hmu(\hf)$ because $(\hf,\hmu)$ is a factor of $(\hF_A,\hm)$.
This proves that
$$
h_\hm(\hF_A,\cW^u) \geq h_\hmu(\hf,W^u).
$$
By Propositions~\ref{p.cuequivalent} and~\ref{p.criterion}, this implies that $\hm$ is
$u$-invariant, as claimed.

\section{Invariant sections and $su$-states}\label{s.states}

We say that an $\hF_A$-invariant probability measure $\hm$ is an \emph{$su$-state} if it is both an $s$-state
and a $u$-state. Here we prove:

\begin{theorem}\label{t.states}
Assume that $A$ admits no invariant section and there exists some periodic point $p$ of $f$ such that
$A^{\per(p)}(p)$ is hyperbolic. Then $A$ has no $su$-states.
\end{theorem}

Assume, by contradiction, that $\hF_A$ does admit some $su$-state $\hm$.
Suppose for a while that $\hm$ admits a \emph{continuous} disintegration $\{\hm_\hx: \hx \in \hM\}$
along the vertical fibers $\{\hx\}\times\pR^2$. The fact that $\hm$ is $\hF_A$-invariant means that
$A(\hx)_*\hm_\hx = \hm_{\hf(\hx)}$ for $\hm$-almost every $\hx$.
Then, by continuity, this must hold for \emph{every} $\hx$.

Let $\hp$ be the fixed point of $\hf$ in $\pi^{-1}(p)$ and $\kappa=\per(p)$ be its period.
Then $\hA^\kappa(\hp)=A^\kappa(p)$ is hyperbolic.
The fact that $\hA^\kappa(\hp)_*\hm_\hp=m_\hp$ implies that $\hm_\hp$ is a convex combination of
not more than two Dirac masses. Then, by $su$-invariance, the same is true about $\hm_\hx$ for every $\hx$.
Thus $\xi(\hx) = \supp\hm_\hx$ defines an invariant section for $\hF_A$, which is in contradiction
with the hypotheses.

In general, disintegrations are only measurable. In what follows we explain how to bypass that and
turn the previous outline into an actual proof of Theorem~\ref{t.states}.

\subsection{Dirac disintegrations}

By the definition of $su$-state, there are disintegrations $\{\hm_\hx^1:\hx\in\hM\}$ and $\{\hm_\hx^2:\hx\in\hM\}$
of $\hm$ and there exists a full $\hmu$-measure subset $U_\hp$ of the neighborhood $V_\hp\approx W^u_\loc(\hp) \times W^s_\loc(\hp)$
such that
\begin{enumerate}
\item[(i)] $(h^ u_{\hx,\hy})_*\hm^1_\hx=\hm^1_\hy$ for every $\hx, \hy \in U_\hp$ with $\hy\in W^u_\loc(\hx)$
($u$-invariance);
\item[(ii)] $\hm^2_\hy=\hm^2_\hz$ for every $\hy, \hz \in V_\hp$ with $\hz\in W^s_\loc(\hy)$
($s$-invariance);
\item[(iii)] $\hm_\hx^1 = \hm_\hx^2$ for every  $\hx \in U_\hp$
(essential uniqueness of disintegrations).
\end{enumerate}
Also, we may choose $U_\hp$ so that $\hm_\hx^1(\Lambda_\hx)=1$
(recall that $\Lambda_\hx=\Lambda \cap (\{\hx\}\times\pR^2)$) for every $\hx\in U_\hp$.

Since the Pesin unstable manifolds $\cW^u(\hz,u)$ vary measurably with the point, we may find compact sets
$\Lambda_1\subset \Lambda_2 \subset \cdots \subset \Lambda$ such that $\hm(\Lambda_j) \to 1$ and
$\cW^u(\hz,u)$ varies continuously on every $\Lambda_j$.
We may choose these compact sets in such a way that $\hF_A(\Lambda_j) \subset \Lambda_{j+1}$ for every $j\ge 1$.
Up to reducing $U_\hp$ if necessary, $\hm^1(\Lambda_{j,\hx}) \to 1$ for every $\hx\in U_\hp$.

Fix any $\hx\in U_\hp$ such that $\hmu^u_\hx(\xi^u(\hx) \setminus U_\hp)=0$. Then define $\hm_\hx = \hm_\hx^1$ and
\begin{enumerate}
\item[(a)] $\hm_\hy = (h^u_{\hx,\hy})_*\hm_\hx$ for every $\hy \in \xi^u(\hx)$;
\item[(b)] $\hm_\hz = \hm_\hy$ for every $\hz \in W^s_\loc(\hy) \cap V_\hp$ with $\hy \in \xi^u(\hx)$.
\end{enumerate}
By (i)-(iii), we have that $\hm_\hy = \hm_\hy^1 = \hm_\hy^2$ for every $\hy \cap \xi^u(\hx) \cap U_\hp$
and $\hm_\hz = \hm_\hz^2$ for every $\hz \in W^s_\loc(\hy) \cap V_\hp$ with $\hy \in \xi^u(\hx) \cap U_\hp$.
By the choice of $\hx$ and the fact that $\hmu$ has local product structure, the latter corresponds to
a full $\hmu$-measure subset of points $\hz \in V_\hp$. In particular, $\{\hm_\hx:\hx\in V_\hp\}$ is a
disintegration of $\hm$ on $V_\hp$.

Let us collect some immediate consequences of the definition of $\hm_\hx$. For $\hx, \hy, \hz$ as in
(a)-(b) above, denote $h^{su}_{\hx,\hz} =  h^s_{\hy,\hz} \circ h^u_{\hx,\hy}$
with $h^s_{\hy,\hz}:\{\hy\}\times\pR^2\to\{\hz\}\times\pR^2$ given by the identity.
For $j\ge 1$, denote $\alpha_j=\hm_\hx(\Lambda_{j,\hx})$; keep in mind that $\alpha_j \to 1$.

\begin{corollary}\label{c.continuity}
For each $j \ge 1$,
\begin{enumerate}
\item[(a)] $\tilde\Lambda_{j,\hz} = h^{su}_{\hx,\hz}(\Lambda_{j,\hx})$ is compact
and varies continuously with $\hz\in V_\hp$;
\item[(b)] the measure $\hm_z \mid \tilde\Lambda_{j,\hz}$ varies continuously with $\hz \in V_\hp$
in the weak$^*$ topology;
\item[(c)] $\hm_\hz(\tilde\Lambda_{j,\hz})=\alpha_j$ for every $\hz \in V_\hp$.
\end{enumerate}
\end{corollary}

Since the matrix $\hA^\kappa(\hp)$ is hyperbolic, its action on the projective space $\pR^2$ is a
North pole-South pole map, that is, a Morse-Smale diffeomorphism with one attractor $a$
and one repeller $r$. We are going to prove:

\begin{proposition}\label{p.pointregidity}
The support of $\hm_\hp$ is contained in $\{a, r\}$.
\end{proposition}

\begin{proof}
Since $\{\hm_\hz:\hz\in V_\hp\}$ is a disintegration and $\hm$ is $\hF_A$-invariant,
\begin{equation}\label{eq.quaseinv}
(\hF_A^\kappa)_*\hm_\hz=\hm_{\hf(\hz)} \text{ for $\hmu$-almost every $\hz\in V_\hp \cap \hf^{-\kappa}(V_\hp)$.}
\end{equation}
The identity may not hold for $\hz=\hp$, but we are going to show that $\hm_\hp$ is at least
``almost $\hF_A$-invariant,'' in a suitable sense:

\begin{lemma}\label{l.invariant}
$\hm_\hp(\hF^{-l \kappa}_A(K)) \geq \hm_\hp(K)$ for any compact set $K \subset \tilde\Lambda_{j,\hp}$
and every $l \geq 1$ and $j\ge 1$.
\end{lemma}

\begin{proof}
Fix $K$, $l$ and $j$. For any $\hq$ close to $\hp$, define
$$
h_\hq = h^s_{\hz,\hF_A^l(\hq)} \circ h^u_{\hy,\hz} \circ h^s_{\hp,\hy}
$$
where $\hy$ and $\hz$ are the points where $W^u_\loc(\hx)$ intersects $W^s_\loc(\hp)$ and
$W^s_\loc(\hf^{l\kappa}(\hq))$, respectively. Keep in mind that the two $s$-holonomies are given
by  the identity. Also, $K \subset\tilde\Lambda_{j,\hp}$ ensures that $h_\hq$ is continuous
restricted to $K$. Define $K_\hq = h_\hq(K)$. Then $K_\hq$ is a compact subset of
$F\{\hf^{l\kappa}(\hq)\}\times\pR^2$ such that $\hm_{\hf^{l\kappa}(\hq)}(K_\hq)=\hm_\hp(K)$.
When $\hq \to \hp$, the point $\hf^{l\kappa}(\hq)$ also goes to $\hp$, and then the same is true
for $\hy$ and $\hz$. Thus $K_\hq \to K$ as $\hq\to\hp$.

Choose $\hq$ close enough to $\hp$ that $\hf^{n\kappa}(\hq) \in V_\hp$ for $0 \le n \le l$
and such that $(\hF_A^{l\kappa})_*\hm_\hq=\hm_{\hf^{l\kappa}(\hq)}$. It follows that
$$
\hm_\hq(\hF^{-l\kappa}_A(K_\hq))
= \hm_{\hf^{l\kappa}(\hq)}(K_\hq) = \hm_\hp(K).
$$
Corollary~\ref{c.continuity}(c) gives that $\hm_\hq(\tilde\Lambda_{k,\hq})=\alpha_k$ for every $k\ge 1$.
Thus
\begin{equation}\label{eq.q}
(\hm_\hq\mid \tilde\Lambda_{k,\hq})(\hF_A^{-l\kappa}(H_\hq))
= \hm_\hq(\tilde\Lambda_{k,\hq}\cap \hF^{-l\kappa}_A(K_\hq))
\ge \hm_\hp(K)+\alpha_k-1.
\end{equation}
By parts (a) and (b) of Corollary~\ref{c.continuity} the compact set $\tilde\Lambda_{k,\hq}$ and
the measure $\hm_\hz \mid \tilde\Lambda_{k,\hq}$ depend continuously on $\hq$.
We know that the same is true for $\hF^{-l\kappa}_A(K_\hq)$.
Thus, making $\hq\to\hp$ in \eqref{eq.q}, we get that
$$
(\hm_\hp \mid \tilde\Lambda_{k,\hp})(\hF^{-l\kappa}_A(K))
\geq \hm_\hp(K)+\alpha_k-1.
$$
Clearly, the left-hand side is less than or equal to $\hm_\hp(\hF^{-l\kappa}_A(K))$.
So, making $k\to\infty$ we get the claim.
\end{proof}

We are ready to complete the proof of Proposition~\ref{p.pointregidity}. Suppose that $\hm_\hp$ is not
supported inside $\{a,r\}$. Then, since the $\tilde\Lambda_{j,\hp}$ are a non-decreasing sequence
whose union has full $\hm_\hp$-measure, for every large $j\ge 1$ the measure $\hm_\hp \mid \tilde\Lambda_{j,\hp}$
is not supported on $\{a, r\})$. Then we can find a compact set $K\subset \tilde\Lambda_{j,\hp}$
contained in a fundamental domain of $\hA^\kappa(\hp)$ with positive $\hm_\hp$-measure.
By Lemma~\ref{l.invariant}, it follows that $\hF_A^{-l\kappa}(K)\ge\hm_\hp(K)>0$ for every $l \geq 0$.
Since these sets are pairwise disjoint, it follows that $\hm_\hp$ is an infinite measure,
which is a contradiction.
\end{proof}

\subsection{Proof of Theorem~\ref{t.states}}

By Proposition~\ref{p.pointregidity}, $\hm_\hp$ is a convex combination of not more than two Dirac
masses. Then, in view of the definition of this disintegration, the same is true about $\hm_\hz$
for every $\hz\in V_\hp$. Then $\hxi(\hz)=\supp\hm_\hz$ defines a continuous map on $V_\hp$ with
values on $\pR^2$ or $\pR^{2,2}$ and such that $\hA(\hz) \hxi(\hz) = \hxi(\hf(\hz))$ for every
$\hz\in V_\hp \cap \hf^{-1}(V_\hp)$.

The same argument shows that for any point $\hy\in\hM$ there exists a continuous disintegration
$\{\hm_{\hy,\hz}:\hz\in V_\hy\}$ of the $su$-state restricted to $V_\hy$. Since disintegrations are
essentially unique and the neighborhoods $V_\hy$ overlap on positive $\hmu$-measure subsets, all
these conditional measures $\hm_{\hy,\hz}$ must be supported on the same number, $1$ or $2$, of
points. Thus, the map $\xi$ in the previous paragraph extends to a continuous invariant section
on the whole $\hM$, which contradicts the assumptions of Theorem~\ref{t.main2}.

\section{Proof of Theorem~\ref{t.main2}}\label{s.main2}

If $\lambda(A)=0$ then, trivially, $A$ is a continuity point. Now assume that $\lambda(A)>0$.
Then (see for instance Kalinin~\cite[Theorem~1.4]{Kal11}) there exists some periodic point $p$
of $f$ such that $A^{\per(p)}(p)$  is hyperbolic. Thus we may use Theorem~\ref{t.states} to
conclude that there are no $su$-states. Now the  proof of continuity of the Lyapunov exponents
is entirely analogous to Section~\ref{s.main1}.

The same arguments also prove the converse: if the  cocycle is hyperbolic at some periodic
point then, again by Theorem~\ref{t.states}, there are no $su$-states and thus the exponent
cannot vanish. The proof of Theorem~\ref{t.main2} is complete.


\appendix

\section{Smooth natural extensions}\label{a.natural}

We show that the natural extension of any $C^k$ local diffeomorphism $f:M\to M$
on a compact manifold admits a $C^k$ realization.

Since $M$ is compact and $f$ is locally injective, we  may find families of open sets
$\{U_i, V_i: i=1, \dots, N\}$ such that: $\{U_1, \dots, U_N\}$ covers $M$;
every $V_i$ contains the closure of $U_i$; and every $f \mid V_i$ is injective.
Take smooth functions $h_i:M\to[0,1]$ such that $h_i \mid U_i \equiv 1$ and
$h_i \mid V_i^c \equiv 0$. Define $h(x)=(h_1(x), \dots, h_N(x))$ for $x\in M$.
Then $h:M\to[0,1]^N$ is such that $h(x) \neq h(y)$ for any pair $(x,y)$ with $x \neq y$
and $f(x) = f(y)$. Since $f$ is locally injective, the set of such pairs is a compact subset
of $M^2$. Hence, there is $\delta>0$ such that $\|h(x) - h(y)\| \ge \delta$ for any $(x,y)$
with $x \neq y$ and $f(x) = f(y)$.

Let $\phi:M\to\RR^m$ be a Whitney embedding of $M$ into some Euclidean space,
and $\psi:M\times D\to\RR^m$ be a tubular neighborhood: $D$ denotes the open unit ball in
$\RR^{m-\dim M}$ and $\psi$ is a smooth  embedding with $\psi(x,0)=\phi(x)$.
Identify $M\times D$ with its image $U=\psi(M\times D)$ through $\psi$.
Fix $\lambda<\delta/4N$ and define
$$
g:M\times D \to M\times D, \quad g(x,v) = (f(x), h(x)/2N+ \lambda v).
$$
It is clear that $g$ is well defined and a $C^k$ local diffeomorphism,
and the image $g(M\times D)$ is relatively compact in $M\times D$.

Suppose that $g(x,v)=g(y,w)$. Then $f(x)=f(y)$ and  $h(x)-h(y) = 2N \lambda (w-v)$.
In particular, $\|h(x)-h(y)\|\le 4 N \lambda < \delta$.
By the definition of $\delta$, this implies that $x=y$.
Then the previous identities imply that $v=w$. This proves that $g$ is injective and,
consequently, an embedding.

For each $\hx=(x_{-n})_n\in\hM$ and $n\ge 1$ the set $g^n(\{x_{-n}\}\times D)$ is a disk
$D_n(\hx)$ of radius $\lambda^n$ inside $\{x_0\}\times D$. These disks are nested and
each $D_{n+1}(\hx)$ is relatively compact in $D_n(\hx)$. Thus, the intersection consists of exactly
one point, which we denote as $\iota(\hx)$. By construction, the map $\iota:\hM\to M\times D$
defined in this way satisfies $g\circ\iota = \iota \circ\hf$. Moreover, the image $\iota(M)$
coincides with $\cap_n g^n(M\times D)$ and so it satisfies $g(\iota(\hM))=\iota(\hM)$.

\end{document}